\documentclass[12pt]{article}
\usepackage{amsmath}
\usepackage{amsfonts}
\usepackage{amssymb}
\usepackage{longtable}

\newtheorem{theorem}{Theorem}
\newtheorem{cor}{Corollary}

\begin{document}

\title{Pseudo-Riemannian manifolds with recurrent spinor fields}

\author{Anton S. Galaev}

\maketitle

\begin{abstract}
The existence of a recurrent spinor field on a pseudo-Riemannian  spin manifold $(M,g)$ is
closely related to the existence of a parallel 1-dimensional complex subbundle of the spinor bundle
 of $(M,g)$. We characterize the following  simply connected
 pseudo-Riemannian manifolds admitting such subbundles in terms of their
 holonomy algebras:  Riemannian manifolds; Lorentzian manifolds;
  pseudo-Riemannian manifolds with irreducible
holonomy algebras; pseudo-Riemannian manifolds of neutral
signature admitting two complementary parallel isotropic
distributions.
\end{abstract}


\section{introduction}

Let $(M,g)$ be a pseudo-Riemannian spin manifold of signature
$(r,s)$ and $S$ the corresponding complex spinor bundle with the
induced connection $\nabla^S$. A spinor field $s\in \Gamma(S)$ is
called {\it recurrent} if
\begin{equation}\label{recs}\nabla^S_Xs=\theta(X)s\end{equation} for all vector fields
$X\in \Gamma(TM)$, here $\theta$ is a complex-valued 1-form. If
$\theta=0$, then $s$ is {\it a parallel} spinor field. For a
recurrent spinor field $s$ there exist locally defined
non-vanishing function $f$ such that $fs$ is parallel  if and only
if $d\theta=0$. If $M$ is simply connected, then such function
exists globally.

The study of Riemannian spin manifolds carrying parallel spinor
fields was initiated by N.\,Hitchin \cite{Hi}, and then by
Th.\,Friedrich \cite{Fr}. M.\,Wang characterized simply connected
Riemannian spin manifolds admitting parallel spinor field in terms
of their holonomy groups \cite{Wang}. A similar result was
obtained by Th.\,Leistner for  Lorentzian manifolds
\cite{L01,LPhD}, by H.\,Baum and I.\,Kath for pseudo-Riemannian
manifolds with irreducible holonomy groups \cite{B-K}, and by
A.\,Ikemakhen in the case of pseudo-Riemannian manifolds of
neutral signature $(n,n)$ admitting two complementary parallel
isotropic distributions \cite{Ik04}.

In \cite{Fr}, Th.\,Friedrich considered Equation \eqref{recs} on a
Riemannian spin manifold with $\theta$ being a real-valued 1-form.
He proved that this equation implies  that the Ricci tensor is
zero and $d\theta=0$. Below we will see that this statement does
not hold for Lorentzian manifolds. Example 1 from \cite{Fr}
provides a solution $s$ to \eqref{recs} with $\theta={\rm
i}\omega$, $d\omega\neq 0$ for a real-valued 1-form $\omega$ on
the compact Riemannian manifold $(M,g)$ being the product of the
non-flat torus $T^2$ and the circle $S^1$. In fact the recurrent
spinor field $s$ comes from a locally defined recurrent spinor
field on the non-Ricci-flat K\"ahler manifold $T^2$ existing by
Theorem~\ref{Th1} below.

The spinor bundle $S$ of $(M,g)$ admits a parallel 1-dimensional
complex subbundle if and only if $(M,g)$ admits non-vanishing
recurrent spinor fields in a neighborhood of each point such that
these fields are proportional on the intersections of the domains
of their definitions. In the present paper we study some classes
of pseudo-Riemannian manifolds $(M,g)$ whose spinor bundles admit
parallel 1-dimensional complex subbundles.

In Section \ref{secR} we prove that if the spinor bundle of a
locally indecomposable Riemannian spin manifold  $(M,g)$ admits a
parallel 1-dimensional complex subbundle, then either $(M,g)$ is a
non-Ricci-flat K\"ahlerian manifold, or it admits a non-zero
parallel spinor field. Moreover, the spinor bundle of any locally
indecomposable non-Ricci-flat K\"ahlerian spin manifold admits
exactly two parallel 1-dimensional complex subbundles. Next, any
simply connected complete  Riemannian spin  manifold with not
irreducible holonomy algebra whose spinor bundle $S$ admits a
parallel 1-dimensional complex subbundle and does not admit any
non-zero parallel spinor field is a product of a non-Ricci flat
K\"ahlerian spin manifold and of a Riemannian spin manifold
admitting a parallel spinor field.

In Section \ref{secL} we prove that  the spinor bundle of an
$n+2$-dimensional simply connected locally indecomposable
Lorentzian spin manifold $(M,g)$ admits a parallel 1-dimensional
complex subbundle if and only if $(M,g)$ admits a parallel
distribution of isotropic lines, i.e. its holonomy algebra
$\mathfrak{g}$ is contained in the parabolic subalgebra
$$\mathfrak{sim}(n)=(\mathbb{R}\oplus\mathfrak{so}(n))\ltimes\mathbb{R}^n\subset\mathfrak{so}(1,n+1),$$
and the subalgebra $\mathfrak{h}={\rm
pr}_{\mathfrak{so}(n)}\mathfrak{g}\subset\mathfrak{so}(n)$
preserves a 1-dimensional complex subspace in the spinor module
$\Delta_n$.

In Section \ref{secpR1} for pseudo-Riemannian spin manifolds with
irreducible holonomy algebras we prove the same statements as for
locally indecomposable Riemannian manifolds in Section \ref{secR}.

In Section \ref{secpR2} we prove that the spinor bundle of any simply connected pseudo-Riemannian spin manifold of neutral
signature $(n,n)$, admitting two complementary parallel isotropic distributions,  admits a parallel 1-dimensional complex subbundle.

In Section \ref{secspinC} we discuss the relation between
recurrent spinor fields and parallel spinor fields of
$spin^C$-bundles that has been studied by A.\,Moroianu
\cite{Mor97} and A.~Ikemakhen \cite{Ik06,Ik07}.

{\it Acknowledgments.} The author is  thankful to Helga Baum for
proposing this problem  and the discussions on this topic.

\section{Preliminaries}
{\bf Spinor representations.} Let us fix some standard denotation.
Let $\mathbb{R}^{r,s}$ be a pseudo-Euclidean space with the metric
$g$ of signature $(r,s)$ ($r$ denotes the number of minuses).  Let
$({\cal C}l_{r,s},\cdot)$ be the corresponding Clifford algebra
and $\mathbb{C} l_{r,s}= {\cal C}l_{r,s}\otimes\mathbb{C}$ be its
complexification. The last algebra can be represented as an matrix
algebra in the following way. Consider the basis
$$\left(u(\epsilon)=\frac{\sqrt{2}}{2}\left(\begin{array}{c}1\\-\epsilon{\rm i}\end{array}\right),\,
\epsilon=\pm 1\right)$$ of $\mathbb{C}^2$. Define the following
isomorphisms of $\mathbb{C}^2$:
$$E={\rm id},\quad T=\left(\begin{array}{cc}0&-{\rm i}\\{\rm i}&0\end{array}\right),
\quad U=\left(\begin{array}{cc}{\rm i}&0\\0&-{\rm
i}\end{array}\right), \quad V=\left(\begin{array}{cc}0&{\rm
i}\\{\rm i}&0\end{array}\right).$$ It holds
$$T^2=-V^2=-U^2=E,\quad UT=-{\rm i} V,\quad VT=i U,\quad UV=-{\rm i} T,$$
$$Tu(\epsilon)=-\epsilon u(\epsilon),\quad Uu(\epsilon)={\rm i}
u(-\epsilon),\quad Tu(\epsilon)=\epsilon u(-\epsilon).$$ Let
$n=r+s$. A basis $e_1,...,e_n$ of $\mathbb{R}^{r,s}$ is called
orthonormal if $g(e_i,e_j)=k_i\delta_{ij}$, where $k_i=-1$ if
$1\leq i\leq r$, and $k_i=1$ if $r+1\leq i\leq n$. Let us fix such
basis. For an integer $m$ denote by $\mathbb{C}(m)$ the algebra of
the complex square matrices of order $m$. Define the following
isomorphisms:
\begin{itemize}
\item[1)] if $n$ is even, then define $\Phi_{r,s}:\mathbb{C} l_{r,s}\to\mathbb{C}\Big(2^{\frac{n}{2}}\Big)$ by
\begin{align}\Phi_{r,s}(e_{2k-1})&=\tau_{2k-1}E\otimes\cdots\otimes E\otimes U\otimes
\underbrace{T\otimes \cdots\otimes T}_{(k-1)-times},\label{Phi1}\\
\Phi_{r,s}(e_{2k})&=\tau_{2k}E\otimes\cdots\otimes E\otimes
V\otimes \underbrace{T\otimes \cdots\otimes
T}_{(k-1)-times},\label{Phi2}\end{align} where $1\leq k\leq
\frac{n}{2}$, $\tau_i={\rm i}$ if $1\leq i\leq r$, and $\tau_i=1$
if $r+1\leq i\leq n$;
\item[2)] if $n$ is odd, then define $\Phi_{r,s}:\mathbb{C} l_{r,s}\to\mathbb{C}\Big(2^{\frac{n-1}{2}}\Big)\oplus\mathbb{C}\Big(2^{\frac{n-1}{2}}\Big)$ by
\begin{align*}\Phi_{r,s}(e_{k})&=(\Phi_{r,s-1}(e_{k}),\Phi_{r,s-1}(e_{k})),\quad k=1,...,n-1,\\
\Phi_{r,s}(e_{n})&=(T\otimes\cdots\otimes T,-{\rm i}
T\otimes\cdots\otimes T).
\end{align*}
\end{itemize}
The obtained representation space
$\Delta_{r,s}=\mathbb{C}^{2^{[\frac{n}{2}]}}$ is called {\it the
spinor module}.
 We write $A\cdot s=\Phi_{r,s}(A)s$ for all
$A\in\mathbb{C} l_{r,s}$, $s\in\Delta_{r,s}$. We will consider the
following basis of  $\Delta_{r,s}$:
$$(u(\epsilon_k,...,\epsilon_1)=u(\epsilon_k)\otimes\cdots\otimes u(\epsilon_1)|\epsilon_i=\pm 1).$$
Recall that the Lie algebra $\mathfrak{spin}(r,s)$ of the Lie
group ${\rm Spin}(r,s)$ can
 be embedded into $\mathbb{C} l_{r,s}$ in the following way:
$$\mathfrak{spin}(r,s)={\rm span}\{e_i\cdot e_j|1\leq i<j\leq n\}.$$
The Lie algebra $\mathfrak{so}(r,s)$ can be identified with the
space of bivectors $\Lambda^2\mathbb{R}^{r,s}$ in such a way that
$$(x\wedge y)z=g(x,z)y-g(y,z)x,\quad x,y,z\in\mathbb{R}^{r,s}.$$
There is the isomorphism
$$\lambda_*:\mathfrak{so}(r,s)\to\mathfrak{spin}(r,s),\quad
\lambda_*(x\wedge y)=x\cdot y.$$ The obtained representation of
$\mathfrak{so}(r,s)$ in $\Delta_{r,s}$ is irreducible if $n$ is
odd, and this representation   splits into the direct some of two
irreducible modules $$\Delta_{r,s}^\pm ={\rm
span}\{u(\epsilon_k,...,\epsilon_1)|\epsilon_1=\cdots=\epsilon_k=\pm
1\}$$ if $n$ is even.

{\bf Holonomy algebras.} Let $(M,g)$ be a pseudo-Riemannian manifold of
 signature $(r,s)$ and $\nabla$ the Levi-Civita connection on
$(M,g)$. Fix a point $x\in M$. The tangent space $T_xM$ can be
identified with the pseudo-Euclidean space $\mathbb{R}^{r,s}$.
Then the holonomy algebra
$\mathfrak{h}_x\subset\mathfrak{so}(T_xM,g_x)$ (i.e. the Lie
algebra of the holonomy group) of $(M,g)$ at the point $x$ is
identified with a subalgebra
$\mathfrak{h}\subset\mathfrak{so}(r,s)$. If $M$ is simply
connected, then the holonomy algebra and the holonomy group
unequally define each other, and we prefer to speak about the
holonomy algebra. If $(M,g)$ is a spin manifold, then it admits a
spinor bundle $S$. The fiber $S_x$ can be identified with the
spinor module $\Delta_{r,s}$. The Levi-Civita connection $\nabla$
defines a connection $\nabla^S$ on the spinor bundle $S$. The
holonomy algebra of this connection coincides with
$\lambda_*(\mathfrak{h})\subset\mathfrak{spin}(r,s)$. The
fundamental property of the holonomy algebra says that if $(M,g)$
is simply connected, then $(M,g)$ admits a parallel vector field
(spinor, distribution, and so on) if and only if $\mathfrak{h}$
annihilates a vector in $\mathbb{R}^{r,s}$ (annihilates a spinor
in $\Delta_{r,s}$, preserves a vector subspace in
$\mathbb{R}^{r,s}$, and so on).

\section{Riemannian manifolds}\label{secR}

\begin{theorem}\label{Th1} Let $(M,g)$ be a locally indecomposable $n$-dimensional simply connected
 Riemannian spin manifold. Then its spinor bundle $S$ admits a parallel 1-dimensional complex subbundle
 if and only if either the holonomy algebra $\mathfrak{h}\subset\mathfrak{so}(n)$ of $(M,g)$
 is one of $\mathfrak{u}(\frac{n}{2})$, $\mathfrak{su}(\frac{n}{2})$,
 $\mathfrak{sp}(\frac{n}{4})$, $G_2\subset\mathfrak{so}(7)$, $\mathfrak{spin}(7)\subset\mathfrak{so}(8)$, or $(M,g)$ is a locally
 symmetric K\"ahlerian manifold. \end{theorem}

{\bf Proof.} Recall that $S$ admits a parallel 1-dimensional
complex subbundle  if and only if the holonomy algebra
$\mathfrak{h}\subset\mathfrak{so}(n)$ preserves a 1-dimensional
complex subspace $l$ of the  spinor module $\Delta_n$. If
$\mathfrak{h}$ annihilates $l$, then $(M,g)$ admits a non-zero
parallel spinor field, and $\mathfrak{h}$ is one of
$\mathfrak{su}(\frac{n}{2})$, $\mathfrak{sp}(\frac{n}{4})$, $G_2$,
$\mathfrak{spin}(7)$ \cite{Wang}. Suppose that $\mathfrak{h}$ does
not annihilate $l$. Since $\mathfrak{h}\subset\mathfrak{so}(n)$,
it is a compact Lie algebra. It is clear that the only compact Lie
algebra admitting a non-trivial exact  2-dimensional real
representation is $\mathfrak{so}(2)$. This shows that
$\mathfrak{h}$ contains  a 1-dimensional center, consequently,
$\mathfrak{h}$ is contained in $\mathfrak{u}(\frac{n}{2})$, i.e.
$(M,g)$ is K\"ahlerian. Since $\mathfrak{su}(\frac{n}{2})$
annihilates two complex 1-dimensional subspaces of $\Delta_n$,
$\mathfrak{u}(\frac{n}{2})$ preserves these two subspaces and it
does not annihilate them. This shows that the spinor bundle $S$ of
a K\"ahlerian manifold admits at least two parallel 1-dimensional
complex subbundles. And it admits exactly two such subbundles if
$\mathfrak{h}=\mathfrak{u}(\frac{n}{2})$. $\Box$

\begin{cor} Let $(M,g)$ be a simply connected  Riemannian spin manifold with irreducible holonomy
algebra and without non-zero parallel spinor fields. Then the
spinor bundle $S$ admits a parallel 1-dimensional complex
subbundle if and only if $(M,g)$ is a  K\"ahlerian  manifold and
it is not Ricci-flat.
\end{cor}

{\it Proof.} The corollary follows from the fact that simply
connected Riemannian manifolds with the holonomy algebras
$\mathfrak{su}(\frac{n}{2}),$
 $\mathfrak{sp}(\frac{n}{4}),$ $G_2$, $\mathfrak{spin}(7)$ do admit non-zero parallel spinor fields.
 $\Box$

\begin{cor} Let $(M,g)$ be a simply connected complete
 Riemannian spin manifold without non-zero parallel spinor
  fields and with not irreducible holonomy algebra.
  Then its spinor bundle $S$ admits a parallel 1-dimensional complex subbundle
if and only if $(M,g)$ is a direct product of a K\"ahlerian not
Ricci-flat spin manifold and of a Riemannian spin manifold with a
non-zero parallel spinor field.
\end{cor}

{\it Proof.} By the de~Rham Theorem, $(M,g)$ can be decomposed as
a direct product of indecomposable simply connected complete
Riemannian manifolds and, probably, of a flat Riemannian manifold,
see e.g. \cite{ESI}. The Riemannian manifolds in this
decomposition are automatically spin, and the spinor bundle $S$ of
$(M,g)$ admits a parallel 1-dimensional complex subbundle if and
only if this is the case for each manifold in the decomposition.
$\Box$

\begin{theorem} Let $(M,g)$ be a locally  indecomposable $n$-dimensional simply
connected non-Ricci-flat K\"ahlerian spin manifold. Then its
spinor bundle $S$ admits exactly two parallel 1-dimensional
complex subbundles. \end{theorem}

This theorem will follow from a more general Theorem \ref{ThpK}
that will be proved  below.

\section{Lorentzian manifolds}\label{secL}

In this section we consider Lorentzian manifolds $(M,g)$, i.e.
pseudo-Riemannian  manifolds of signature $(1,n+1)$, $n\geq 0$.
Holonomy algebras of Lorentzian spin manifold admitting non-zero
parallel spinor fields are classified in \cite{L01,LPhD}. We
suppose now that the spinor bundle of $(M,g)$ admits a parallel
1-dimensional complex subbundle and $(M,g)$ does not admit any
parallel spinor.

\begin{theorem}\label{decLor}
Let $(M,g)$ be a simply connected complete Lorentzian spin
manifold. Suppose that $(M,g)$ does not admit a parallel spinor.
Then its spinor bundle $S$ admits a parallel 1-dimensional complex
subbundle if and only if one of the following conditions holds:

\begin{itemize}
\item[1)] $(M,g)$ is a direct product of $(\mathbb{R},-(dt)^2)$ and of a Riemannian spin
manifold $(N,h)$ such that the spinor bundle of $(N,h)$  admits a
parallel 1-dimensional complex subbundle and $(N,h)$ does not
admit any non-zero parallel spinor field;

\item[2)]
 $(M,g)$ is a direct product of an
indecomposable Lorentzian spin manifold and of Riemannian spin
manifold $(N,h)$ such that the spinor bundles of both manifolds
admit parallel 1-dimensional complex subbundles and at least one
of these manifolds does not admit any non-zero parallel spinor
field.
\end{itemize} \end{theorem}

{\bf Proof.} The proof of this theorem follows directly from the
Wu decomposition Theorem, the proof of a similar statement can be
found in \cite{L01,LPhD}. $\Box$

Note that if $(M,g)$ admits a parallel spinor field, then it
satisfies the analog of Theorem \ref{decLor} with $(N,h)$
admitting a parallel spinor field \cite{L01,LPhD}.

Now we may consider locally indecomposable Lorentzian manifolds
$(M,g)$. Suppose that the spinor bundle of $(M,g)$  admits a
parallel 1-dimensional complex subbundle $l$. Let  $s\in
\Gamma(l)$ be a local non-vanishing section of $l$.
 Let $p\in\Gamma(TM)$ be its Dirac current. Recall that $p$ is defined from the equality
 $$g(p,X)=-<X\cdot s,s>,$$ where $<,>$ is a Hermition product on $S$.
We claim that $p$ is a recurrent vector field. Indeed, following the computations from \cite{LPhD},
for any vector fields $X$ and $Y$ we get \begin{align*}g(\nabla_Yp,X)&=Y(g(p,X))-g(p,\nabla_YX)=-(Y(<X\cdot s,s>)-<\nabla_YX\cdot s,s>)\\
&=-(<\nabla_YX\cdot s,s>+<X\cdot\nabla^S_Ys,s>+<X\cdot s,\nabla^S_Ys>-<\nabla_YX\cdot s,s>)\\
&= -(\theta(Y)+\overline{\theta(Y)})<X\cdot
s,s>=(\theta(Y)+\overline{\theta(Y)})g(p,X).\end{align*} We obtain
\begin{equation}\label{nablap}\nabla_Yp=2{\rm Re}(\theta)(Y)p.\end{equation} Recall that in the
Lorentzian signature the Dirac current satisfies $g(p,p)\leq 0$
and the zeros of $p$ coincide with the zeros of $s$. Since $s$ is
non-vanishing and $p$ is recurrent, we see that either $g(p,p)<0$,
or $g(p,p)=0$. In the first case the manifold is decomposable.
Thus we get that $p$ is an isotropic recurrent vector field. We
conclude that $(M,g)$ admits a parallel distribution of isotropic
lines. The holonomy algebras of such manifolds are classified
\cite{BB-I,L07,G06,ESI}. We recall this classification now.

Let $(M,g)$ be a locally indecomposable Lorentzian manifold of
dimension $n+2$. Suppose that $(M,g)$ admits a parallel
distribution of isotropic lines. The tangent space to $(M,g)$ at
any point can be identified with the
 Minkowski space $\mathbb{R}^{1,n+1}$.  We fix a basis $p,e_1,...,e_n,q$ of $\mathbb{R}^{1,n+1}$ such that the only non-zero values of the metric are
$g(p,q)=g(q,p)=1$ and $g(e_i,e_i)=1$.  We will denote by
$\mathbb{R}^n\subset\mathbb{R}^{1,n+1}$ the Euclidean subspace
spanned by the vectors $e_1,...,e_n$. Suppose that $\mathbb{R} p$
corresponds to the parallel distribution of isotropic lines.
Denote  by $\mathfrak{sim}(n)$ the subalgebra of
$\mathfrak{so}(1,n+1)$ that preserves the isotropic line
$\mathbb{R} p$. The Lie algebra $\mathfrak{sim}(n)$ can be
identified with the following matrix algebra:
$$\mathfrak{sim}(n)=\left\{\left. \left (\begin{array}{ccc} a
&X^t & 0\\ 0 & A &-X \\ 0 & 0 & -a \\
\end{array}\right)\right|\, a\in \mathbb{R},\,A \in \mathfrak{so}(n),\, X\in \mathbb{R}^n \right\}
.$$ The above matrix can be identified with the triple $(a,A,X)$.
We get the decomposition
$$\mathfrak{sim}(n)=(\mathbb{R}\oplus\mathfrak{so}(n))\ltimes\mathbb{R}^n,$$ which means that $\mathbb{R}\oplus\mathfrak{so}(n)\subset\mathfrak{sim}(n)$ is a subalgebra and
$\mathbb{R}^n\subset\mathfrak{sim}(n)$ is an ideal, and the Lie
brackets of $\mathbb{R}\oplus\mathfrak{so}(n)$ with $\mathbb{R}^n$
are given by the standard representation of
$\mathbb{R}\oplus\mathfrak{so}(n)$ in $\mathbb{R}^n$. The holonomy
algebra $\mathfrak{g}$ of $(M,g)$ is contained in
$\mathfrak{sim}(n)$. Any subalgebra
$\mathfrak{h}\subset\mathfrak{so}(n)$ can be decomposed into the
direct sum
$\mathfrak{h}=\mathfrak{h}'\oplus\mathfrak{z}(\mathfrak{h})$,
where $\mathfrak{h}'=[\mathfrak{h},\mathfrak{h}]$ and
$\mathfrak{z}(\mathfrak{h})$ is the center of $\mathfrak{h}$.

\begin{theorem}\label{holclas}
A subalgebra $\mathfrak{g}\subset \mathfrak{sim}(n)$ is the
holonomy algebra of a locally indecomposable Lorentzian manifold
if and only if it is conjugated to one of the following
subalgebras:
\begin{description}
\item[type 1.] $\mathfrak{g}^{1,\mathfrak{h}}=(\mathbb{R}\oplus\mathfrak{h})\ltimes\mathbb{R}^n$, where $\mathfrak{h}\subset\mathfrak{so}(n)$ is the holonomy algebra of a Riemannian manifold;

\item[type 2.] $\mathfrak{g}^{2,\mathfrak{h}}=\mathfrak{h}\ltimes\mathbb{R}^n,$ where $\mathfrak{h}\subset\mathfrak{so}(n)$ is the holonomy algebra of a Riemannian manifold;

\item[type 3.] $\mathfrak{g}^{3,\mathfrak{h},\varphi}=\{(\varphi(A),A,0)|A\in\mathfrak{h}\}\ltimes\mathbb{R}^n,$ where $\mathfrak{h}\subset\mathfrak{so}(n)$ is
the holonomy algebra of a Riemannian manifold
 with $\mathfrak{z}(\mathfrak{h})\neq\{0\}$, and  $\varphi :\mathfrak{h}\to\mathbb{R}$ is a non-zero linear map with
$\varphi|_{\mathfrak{h}'}=0$;

\item[type 4.] $\mathfrak{g}^{4,\mathfrak{h},m,\psi}=\{(0,A,X+\psi(A))|A\in\mathfrak{h},X\in\mathbb{R}^m\},$ where  $0<m<n$ is an integer,
$\mathfrak{h}\subset\mathfrak{so}(m)$ is the holonomy algebra of a
Riemannian manifold with $\dim\mathfrak{z}(\mathfrak{h})\geq n-m$,
and $\psi:\mathfrak{h}\to \mathbb{R}^{n-m}$ is a surjective linear
map with $\psi|_{\mathfrak{h}'}=0$.
\end{description}\end{theorem}

Recall that for the subalgebra
$\mathfrak{h}\subset\mathfrak{so}(n)$
 there is an orthogonal decomposition
$$\mathbb{R}^{n}=\mathbb{R}^{n_1}\oplus\cdots\oplus\mathbb{R}^{n_s}\oplus\mathbb{R}^{n_{s+1}}$$
and the corresponding decomposition into the direct sum of ideals
\begin{equation}\label{LM0B}\mathfrak{h}=\mathfrak{h}_1\oplus\cdots\oplus\mathfrak{h}_s\oplus\{0\}\end{equation} such that $\mathfrak{h}$ annihilates
$\mathbb{R}^{n_{s+1}}$, $\mathfrak{h}_i(\mathbb{R}^{n_j})=0$ for
$i\neq j$, and $\mathfrak{h}_i\subset\mathfrak{so}(n_i)$ is an
irreducible subalgebra for $1\leq i\leq s$. Moreover, the
subalgebras $\mathfrak{h}_i\subset\mathfrak{so}(n_i)$ are the
holonomy algebras of Riemannian manifolds.

In \cite{L01,LPhD} it is proved that $(M,g)$ admits a non-zero
parallel spinor field if and only if
$\mathfrak{g}=\mathfrak{g}^{2,\mathfrak{h}}=\mathfrak{h}\ltimes\mathbb{R}^n$
and in the decomposition \eqref{LM0B} of
$\mathfrak{h}\subset\mathfrak{so}(n)$ each subalgebra
$\mathfrak{h}_i\subset\mathfrak{so}(n_i)$ coincides with one of
the Lie algebras
 $\mathfrak{su}(\frac{n_i}{2})$, $\mathfrak{sp}(\frac{n_i}{4})$, $G_2\subset\mathfrak{so}(7)$, $\mathfrak{spin}(7)\subset\mathfrak{so}(8)$.

\begin{theorem}\label{Threc5} Let $(M,g)$ be a simply connected locally
indecomposable $(n+2)$-dimensional Lorentzian spin manifold. Then
its spinor bundle $S$ admits a parallel 1-dimensional complex
subbundle if and only if $(M,g)$ admits a parallel distribution of
isotropic lines (i.e. its holonomy algebra $\mathfrak{g}$ is
contained in $\mathfrak{sim}(n)$), and  in the decomposition
\eqref{LM0B} of $\mathfrak{h}={\rm
pr}_{\mathfrak{so}(n)}\mathfrak{g}$ each subalgebra
$\mathfrak{h}_i\subset\mathfrak{so}(n_i)$ coincides with one of
the Lie algebras
  $\mathfrak{u}(\frac{n_i}{2}),\mathfrak{su}(\frac{n_i}{2}),
 \mathfrak{sp}(\frac{n_i}{4}),G_2,\mathfrak{spin}(7)$,
or with the holonomy algebra of an indecomposable K\"ahlerian
symmetric space. The number of parallel 1-dimensional complex
subbundles of $S$ equals to the number of $1$-dimensional complex
subspaces of $\Delta_n$ preserved by $\mathfrak{h}$.
\end{theorem}

{\bf Proof.} It is enough to find for each Lie algebra
$\mathfrak{g}$ from Theorem \ref{holclas} all  1-dimensional
complex invariant submodules of $\Delta_{1,n+1}$.

Under the identification $\mathfrak{so}(1,n+1)\simeq
\Lambda^2\mathbb{R}^{1,n+1}$, an element
$(a,A,X)\in\mathfrak{sim}(n)$ corresponds to the bivector
$$-ap\wedge q+A-p\wedge X.$$ Recall that $$\Delta_{1,n+1}\simeq
\Delta_{n} \otimes \Delta_{1,1},\quad \Delta_{1,1}\simeq
\mathbb{C}^2.$$ Consider the vectors $e_-=\frac{\sqrt{2}}{2}(p-q)$
and $e_+=\frac{\sqrt{2}}{2}(p+q)$ and the orthonormal basis $e_-
,e_+,e_1,...,e_n$ of $\mathbb{R}^{1,n+1}$. Note that
$p=\frac{\sqrt{2}}{2}(e_-+e_+)$ and
$q=\frac{\sqrt{2}}{2}(e_--e_+)$.

Let us find all complex 1-dimensional
$p\wedge\mathbb{R}^n$-invariant submodules of $\Delta_{1,n+1}$.
Let $w\in \Delta_{1,n+1}$ be a spinor such that $(p\wedge
e_i)\cdot w=c_i w,$ $c_i\in\mathbb{C}$, for all $i=1,...,n$. We
may write
$$w=w_+\otimes u(1)+w_-\otimes u(-1),$$ where
$w_\pm\in\Delta_{n}$. Using \eqref{Phi1}, \eqref{Phi2} and the
computations from \cite{L01,LPhD} it is easy to get that $$
(e_i\wedge p)\cdot w=\frac{\sqrt{2}}{2}e_i\cdot(e_-+e_+) \cdot
w=\sqrt{2}(e_i\cdot w_-)\otimes u(1).$$ Hence the equalities
$(p\wedge e_i)\cdot w=c_i w$ imply $c_i=0$, $e_i\cdot w_-=0$.
Since the vectors $e_i$ act in $\Delta_n$ as isomorphisms, we get
$w_-=0$. Thus, $w=w_+\otimes u(1)$, and all complex 1-dimensional
$p\wedge\mathbb{R}^n$-invariant submodules of $\Delta_{1,n+1}$ are
contained in $\Delta_{n}\otimes u(1)$. Moreover, this is a trivial
$p\wedge\mathbb{R}^n$-module. The same statement holds for
$p\wedge\mathbb{R}^m$ in the case of the Lie algebra
$\mathfrak{g}^{4,\mathfrak{h},m,\psi}$. Next, $$(p\wedge
q)\cdot(w_+\otimes u(1))=2w_+\otimes u(1),\quad A(w_+\otimes
u(1))=A(w_+)\otimes u(1)$$ for all $w_+\in\Delta_n$ and
$A\in\mathfrak{so}(n)$. This shows that all
$\mathfrak{g}$-invariant $1$-dimensional subspaces  of
$\Delta_{1,n+1}$ are of the form $l\otimes u(1)$, where $l\subset
\Delta_n$ is an $\mathfrak{h}$-invariant $1$-dimensional subspace.
The theorem is true. $\Box$

Examples of analytical Lorentzian manifolds with each possible
holonomy algebra, in particular, with the holonomy algebras as in
Theorem \ref{Threc5} can be obtained using the construction from
\cite{G06}.

Let us consider one example. Let $(M,g)$ be a simply connected
Lorentzian spin manifold with the holonomy algebra
$\mathfrak{g}=(\mathbb{R}\oplus\mathfrak{h})\ltimes\mathbb{R}^n$,
where $\mathfrak{h}\subset\mathfrak{so}(n)$ is the holonomy
algebra of a Riemannian manifold carrying a non-zero parallel
spinor field. By Theorem \ref{Threc5}, on $(M,g)$ locally exists a
recurrent spinor field $s$ (restricting the consideration to a
local subset of $(M,g)$ with the same holonomy algebra, we may
assume that $s$ is defined globally). By the assumption on
$\mathfrak{h}$ and by the proof of Theorem \ref{Threc5}, the
corresponding 1-form $\theta$ must be real-valued. By the fact
that the manifolds with the holonomy algebra $\mathfrak{g}$ do not
admit any non-zero parallel vector field, and by \eqref{nablap},
$d\theta\neq 0$. This shows that the statement of Theorem 1 from
\cite{Fr} mentioned in Introduction does not hold for Lorentzian
manifolds.

\section{Pseudo-Riemannian manifolds with irreducible holonomy
algebras}\label{secpR1}

\begin{theorem}\label{Th1pR} Let $(M,g)$ be a  simply connected
 pseudo-Riemannian spin manifold of signature $(r,s)$  with an irreducible holonomy algebra $\mathfrak{h}\subset\mathfrak{so}(r,s)$.
 Then its spinor bundle $S$ admits a parallel 1-dimensional complex subbundle
 if and only if either the holonomy algebra $\mathfrak{h}$
 is one of $\mathfrak{u}(\frac{r}{2},\frac{s}{2})$, $\mathfrak{su}(\frac{r}{2},\frac{s}{2})$,
 $\mathfrak{sp}(\frac{r}{4},\frac{s}{2})$, $G_2\subset\mathfrak{so}(7)$, $G^*_{2(2)}\subset\mathfrak{so}(3,4)$, $G^\mathbb{C}_2\subset\mathfrak{so}(7,7)$,
 $\mathfrak{spin}(7)\subset\mathfrak{so}(8)$, $\mathfrak{spin}(3,4)\subset\mathfrak{so}(4,4)$, $\mathfrak{spin}(7,\mathbb{C})\subset\mathfrak{so}(8,8)$,  or $(M,g)$ is a locally
 symmetric pseudo-K\"ahlerian manifold. \end{theorem}

{\bf Proof.} Recall that irreducible holonomy algebras of simply
connected pseudo-Riemannian manifolds $(M,g)$ admitting non-zero
parallel spinor fields are exhausted by
$\mathfrak{su}(\frac{r}{2},\frac{s}{2})$,
 $\mathfrak{sp}(\frac{r}{4},\frac{s}{2})$, $G_2\subset\mathfrak{so}(7)$, $G^*_{2(2)}\subset\mathfrak{so}(3,4)$, $G^\mathbb{C}_2\subset\mathfrak{so}(7,7)$,
 $\mathfrak{spin}(7)\subset\mathfrak{so}(8)$, $\mathfrak{spin}(3,4)\subset\mathfrak{so}(4,4)$,
 $\mathfrak{spin}(7,\mathbb{C})\subset\mathfrak{so}(8,8)$ \cite{B-K}. From \cite{Armstrong} it  follows that manifolds with each of these holonomy algebras are Ricci-flat.

 Suppose that $(M,g)$ is not Ricci-flat and its spinor bundle $S$ admits a parallel 1-dimensional complex
 subbundle $l$. Let $s$ be a local section of $l$ defined at a
 point $y\in M$. Equation \eqref{recs} easily implies
 \begin{equation}\label{Rs} R_y^S(X,Y)s_y=(d\theta)_y(X,Y)s_y\end{equation} for
 all  $X,Y\in T_yM$. Here $R^S$ is the curvature tensor of the
 connection $\nabla^S$. This equation shows that $R_y^S(X,Y)$ acts
 on $l_y$ as a multiplications by complex numbers. This and the
 Ambrose-Singer Theorem imply that  the
 holonomy algebra $\mathfrak{h}_x$, $x\in M$, of the manifold $(M,g)$ acts on $l_x$  by multiplications by
  complex numbers and
 this action is non-trivial.  Since $\mathfrak{h}\subset\mathfrak{so}(r,s)$ is irreducible,
  it is a reductive Lie
 algebra. We get that $\mathfrak{h}$
  contains a commutative ideal. By the Schur Lemma, this ideal is generated a the complex
 structure  $J\in\mathfrak{u}(\frac{r}{2},\frac{s}{2})$.
  Consequently, $\mathfrak{h}\subset\mathfrak{u}(\frac{r}{2},\frac{s}{2})$. Since
 $\mathfrak{su}(\frac{r}{2},\frac{s}{2})$ annihilates two 1-dimensional
 complex subspaces in $\Delta_{r,s}$ \cite{B-K},
 $\mathfrak{u}(\frac{r}{2},\frac{s}{2})$ preserves these subspaces.
 Consequently the spinor bundle of each pseudo-K\"ahlerian
 manifold  admits at least two parallel 1-dimensional complex
 subbundles. $\Box$

\begin{cor} Let $(M,g)$ be a simply connected  pseudo-Riemannian spin manifold with irreducible holonomy
algebra and without non-zero parallel spinor fields. Then the
spinor bundle $S$ admits a parallel 1-dimensional complex
subbundle if and only if $(M,g)$ is a  pseudo-K\"ahlerian manifold
and it is not Ricci-flat.
\end{cor}

\begin{theorem}\label{ThpK} Let $(M,g)$ be a simply
connected not Ricci-flat pseudo-K\"ahlerian spin manifold with an
irreducible holonomy algebra. Then its spinor bundle $S$ admits
exactly two parallel 1-dimensional complex subbundles.
\end{theorem}

{\bf Proof.} Let $(M,g)$ be a simply connected not Ricci-flat
pseudo-Riemannian spin manifold with an irreducible holonomy
algebra. Suppose that the spinor bundle $S$ admits a parallel
1-dimensional complex subbundle $l$. Let $s$ be a local
non-vanishing section of $l$, then $s$ is a recurrent spinor
field. In the proof of Theorem \ref{Th1pR} we have seen that
$\mathfrak{h}$ contains a
 1-dimensional commutative ideal generated by the complex
 structure  $J\in\mathfrak{u}(\frac{r}{2},\frac{s}{2})$ and in Equality \eqref{Rs} $R_x(X,Y)$ acts
 on $s_x$ as its projection on $\mathbb{R} J\subset\mathfrak{h}_x$. Consequently Equality
 \eqref{Rs} implies that $(d\theta(X,Y))^2=-a(X,Y)^2$, $a(X,Y)\in\mathbb{R}$, that is
 $d\theta=\i \omega$ for a real-valued 2-form $\omega$ on $M$.

 Now we follow the  ideas form \cite{Mor97} and \cite{Ik06}.
  The computations similar to
 \cite[Lem. 3.1]{Mor97} and \cite[Th. 2]{Fr} show the equality \begin{equation}\label{M2}
 {\rm Ric}(X)\cdot s={\rm i}(X\lrcorner \omega)\cdot s\end{equation}
 for all vector fields $X$. Here ${\rm Ric}$ is the  Ricci operator.
Consider the  distributions $T$ and $E$ defined in the following
way: $$T_x=\{X\in T_xM|X\cdot s_x=0\},$$ $$E_x=\{X\in T_xM|\exists
Y\in T_xM,\,X\cdot s_x={\rm i} Y\cdot s_x\}.$$ Since we may choose
$s$ in a neighborhood of each point and any two such spinor fields
are proportional on the intersections of the domains of their
definitions, the distributions $T$ and $E$ are defined globally on
$M$. We claim that the both distributions are parallel. Suppose
that $X\in \Gamma(T)$, then for any local vector field $Z$ it
holds $$(\nabla_ZX)\cdot s=\nabla^S_Z(X\cdot
s)-X\cdot(\nabla^S_Zs)=0-X\cdot(\theta(Z)s)=-\theta(Z)(X\cdot
s)=0,$$ i.e. $T$ is parallel. Let $X\in \Gamma(E)$, then locally
there exist vector field $Y$ such that $X\cdot s={\rm i} Y\cdot
s$. For any local vector field $Z$ it holds $$(\nabla_ZX)\cdot
s=\nabla^S_Z(X\cdot s)-X\cdot(\nabla^S_Zs)=\nabla^S_Z({\rm i}
Y\cdot s)-X\cdot(\theta(Z)s)={\rm i} (\nabla_Z
Y+\theta(Z)Y-\theta(Z)Y)\cdot s={\rm i} (\nabla_Z Y)\cdot s,$$
i.e. $E$ is parallel. Since ${\rm Ric}\neq 0$ and $w=id\theta\neq
0$, \eqref{M2} shows that $T\neq TM$ and $E\neq 0$. Since the
holonomy algebra of $(M,g)$ is irreducible, $T=0$ and $E=TM$. Thus
for any local vector field  $X$ there exists a unique vector field
$Y$ such that $X\cdot s={\rm i} Y\cdot s$. This defines an
endomorphism $I$ of the tangent bundle such that
\begin{equation}\label{M7} X\cdot s={\rm i} I(X)\cdot s.\end{equation}
Consequently, $$X\cdot s={\rm i} I(X)\cdot s=- I^2(X)\cdot s,\quad
(X+I^2(X))\cdot s=0.$$ Since $T=0$, it holds $I^2(X)=-X$, i.e. $I$
is an almost complex structure on $M$. Now we show that $I$ is
$g$-orthogonal. It holds \begin{multline*}2g(I(X),I(Y))s=I(X)\cdot
I(Y)\cdot s+I(Y)\cdot I(X)\cdot s=-{\rm i}(I(X)\cdot Y\cdot
s+I(Y)\cdot X\cdot s)\\=-{\rm i}(2g(I(X),Y)s-Y\cdot I(X)\cdot
s+2g(X,I(Y))s-X\cdot I(Y)\cdot s)\\=X\cdot Y\cdot s+Y\cdot X\cdot
s-2{\rm i}(g(I(X),Y)+g(X,I(Y)))s\\=2g(X,Y)s-2{\rm
i}(g(I(X),Y)+g(X,I(Y)))s.\end{multline*} This implies, in
particular, $g(I(X),I(Y))=g(X,Y)$. Next we claim that $I$ is
parallel. Indeed, \eqref{M7} implies for all vector fields $X$ and
$Y$ the following:
$$\nabla_YX\cdot s+X\cdot \theta(Y)\cdot s={\rm i}\nabla_Y(IX)\cdot s+{\rm i} I(X)\cdot\theta(Y)\cdot s.$$
Hence, $$\theta(Y) X\cdot  s={\rm i}\nabla_Y(IX)\cdot s.$$ Using
this and \eqref{M7} applied to $\nabla_Y X$, we get
$$0={\rm i}(\nabla_Y(IX)-I(\nabla_YX))\cdot s.$$ This shows that $\nabla_Y(IX)-I(\nabla_YX)\in\Gamma(T)$. Hence, $\nabla_Y(IX)-I(\nabla_YX)=0$, i.e.  $\nabla I=0$.
Thus, $s$ defines a K\"ahler structure $I$. Note that this gives
another proof of Theorem \ref{Th1pR}.

Suppose that we have two recurrent spinor fields $s$ and $s_1$.
These fields define two K\"ahler  structures $I$ and $I_1$
satisfying \eqref{M7} for $s$ and $s_1$, respectively. Since $I$
and $I_1$ are parallel, their values $I_x$ and $I_{1x}$ at a point
$x\in M$ commute with the holonomy algebra
$\mathfrak{h}_x\subset\mathfrak{so}(T_xM,g_x)\simeq\mathfrak{so}(r,s)$.
Since $\mathfrak{h}_x\subset\mathfrak{so}(T_xM,g_x)$ is
irreducible, by the Schur Lemma, the centralizer of
$\mathfrak{h}_x$ in $\mathfrak{so}(T_xM,g_x)$ is isomorphic ether
to ${\rm i}\mathbb{R}$, or to $\mathfrak{sp}(1)$. In the last
case, $\mathfrak{h}_x$ must be contained in
$\mathfrak{sp}(\frac{r}{4},\frac{s}{4})$. This gives a
contradiction, since the manifolds with such holonomy algebra are
Ricci-flat. We conclude that  $I_x=\alpha I_{1x}$ for some
$\alpha\in\mathbb{R}$. Since $I_x^2=I^2_{1x}=-{\rm id}$,
$\alpha=\pm 1$. Since the both tensor fields $I$ and $I_1$ are
parallel, $I=\alpha I_1$.

Suppose that $I=I_1$. Let $\Omega$ be the pseudo-K\"ahlerian form
corresponding to $I$. Let $e_1,...,e_m,I e_1,...,I e_m$
($r+s=n=2m$) be a local orthonormal basis of $\Gamma(TM)$, i.e.
$g(e_i,e_i)=g(Ie_i,Ie_i)=k_i$, where $k_i=-1$ for
$i=1,...,\frac{r}{2}$, and $k_i=1$ for $i=\frac{r}{2}+1,...,m$.
Let $e^1,...,e^m,(I e)^1,...,(I e)^m$ be the dual basis of
$\Gamma(TM)$. Then $$\Omega=-\sum_{i=1}^m k_i(I e)^i\wedge e^i.$$
Using this and \eqref{M7}, we get
$$\Omega\cdot s=-\sum_{i=1}^m k_i (I e_i)\cdot e_i\cdot s=-{\rm i}
\sum_{i=1}^m k_i(I e_i)\cdot (Ie_i)\cdot s=-{\rm i} m s.$$
Similarly, $\Omega\cdot s_1=-{\rm i} m s_1$. Recall that the
spinor bundle $S$ of a pseudo-K\"ahlerian manifold admits the
decomposition
$$S=\sum_{k=0}^m S_k,\quad (S_k)_x\simeq\Lambda^k \mathbb{C}^m,$$ where $S_k$ is the
eigenspace corresponding to the eigenvalue $(m-2k) {\rm i}$ of the
operator of the Clifford multiplication by the form $\Omega$. We
get that $s,s_1\in S_0$. Since $S_0$ is of rank 1, $s$ and $s_1$
are proportional and they belong to the same 1-dimensional
parallel subbundle of $S$. Above we have seen that $S$  admits at
least two parallel 1-dimensional complex subbundles. We see now
that it admits exactly two subbundles and the second one
corresponds to the K\"ahler structure $-I$ (and coincides with
$S_m$ that has dimension $1$). $\Box$

\section{Pseudo-Riemannian manifolds of neutral signature}\label{secpR2}

\begin{theorem} Let $(M,g)$ be a simply connected pseudo-Riemannian spin manifolds of neutral
signature $(n,n)$ admitting two complementary parallel isotropic distributions. Then
 the spinor bundle $S$ of $(M,g)$ admits a parallel 1-dimensional complex subbundle.
\end{theorem}

{\bf Proof.} The tangent space of such manifold can be identified
with the space $\mathbb{R}^{n,n}$ and it admits the decomposition
$\mathbb{R}^{n,n}=W\oplus W_1$, where $W,W_1\subset
\mathbb{R}^{n,n}$ are isotropic. The holonomy algebra of such
manifold is of the form
$$
\widetilde{\mathfrak{h}}
=\left\{A=\left.\left(\begin{array}{cc}B&0\\0&-B^T\end{array}\right)\right|B\in\mathfrak{h}\right\}$$
for a subalgebra $\mathfrak{h}\subset\mathfrak{gl}(n,\mathbb{R})$,
which is the restriction of $\widetilde{\mathfrak{h}}$ to $W$. In
\cite{Ik04} it is shown that $\widetilde{\mathfrak{h}}$
annihilates an element in $\Delta_{n,n}$ if and only if
$\mathfrak{h}\subset\mathfrak{sl}(n,\mathbb{R})$. It is also
proved that the above element $A\in\widetilde{\mathfrak{h}}$ acts
in $\Delta_{n,n}$ as
$$\frac{1}{2}n\cdot +\frac{1}{4}\sum_{i=1}^n(e_i^*\cdot A(e_i)-A(e_i^*)\cdot e_i)\cdot,$$
where $e_1,...,e_n$ and $e_1^*,...,e_n^*$ are bases of $W$ and
$W_1$, respectively, such that $g(e_i,e^*_j)=\delta_{ij}$. Hence
$A_0=\left(\begin{array}{cc}E&0\\0&-E\end{array}\right)$ acts in
$\Delta_{n,n}$ as the multiplication by $n$. Since the Lie algebra
$\widetilde{\mathfrak{sl}(n,\mathbb{R})}$ annihilates some
non-zero elements of $\Delta_{n,n}$, the Lie algebra
$\widetilde{\mathfrak{gl}(n,\mathbb{R})}=\widetilde{\mathfrak{sl}(n,\mathbb{R})}\oplus\mathbb{R}
A_0$ preserves the corresponding lines in $\Delta_{n,n}$. This
proves the theorem. $\Box$

\section{Relation to parallel spinor fields of $spin^C$ bundles}\label{secspinC}

Let $(M,g)$ be a simply connected  pseudo-Riemannian spin manifold
with the holonomy algebra $\mathfrak{h}\subset\mathfrak{so}(p,q)$
and with the spin bundle $S$. Let $S^C$ be the $spin^C$ bundle on
$(M,g)$ given by a complex line bundle $L$ over $(M,g)$ with a
connection 1-form ${\rm i} A$ and with the holonomy algebra
$\mathfrak{h}_L\subset{\rm i} \mathbb{R}$. The holonomy algebra
$\widehat{\mathfrak{h}}$ of the induced connection on $S^C$ is
contained in $\lambda_*(\mathfrak{h})\oplus {\rm
i}\mathbb{R}\subset \mathfrak{spin}(p,q)\oplus{\rm i} \mathbb{R}$.
In \cite{Ik06} it is proved that the existence of a non-zero
parallel spinor field in $S^C$ is equivalent to the existence of a
spinor $v\in\Delta_{p,q}$ such that $\xi\cdot v={\rm i} t v$ for
all $(\xi,{\rm i} t)\in\widehat{\mathfrak{h}}$. Hence, the
existence of a parallel spinor field in $S^C$ implies the
existence of a parallel 1-dimensional complex subbundle of $S$. In
general the converse statement is not true. For example, in the
Lorentzian case the existence of a spinor field in $S^C$ implies
the existence of a parallel  vector field on $(M,g)$ \cite{Ik07},
while we have seen above that there exist Lorentzian manifolds
admitting recurrent spinor fields and no non-zero parallel vector
fields.

\end{document}